\documentstyle{amsart}

\def\ds{\displaystyle}

\newsymbol\Subset 1362
\def\ds{\displaystyle}

\newcommand{\R}{{\Bbb R}}

\newcommand{\M}{{\Bbb M}}

\def\XXint#1#2#3{{\setbox0=\hbox{$#1{#2#3}{\int}$ } 
\vcenter{\hbox{$#2#3$ }}\kern-.6\wd0}}

\def\ol{\overline}

\def\cA{{\mathcal A}}
\def\B{{\mathcal B}}

\def\bA{{\mathbf A}}

\def\div{{\rm div}}

\def\ba{{\mathbf a}}
\def\bu{{\mathbf u}}
\def\bz{{\mathbf z}}

\def\ff{{\mathbf f}}
\def\bb{{\mathbf b}}
\def\bv{{\mathbf v}}
\newcommand{\Mb}{{\mathbb M}}
\def\diam{{\rm diam\,}}
\def\sign{{\rm sign\,}}

\newtheorem{thm}{Theorem}
\newtheorem{lem}[thm]{Lemma}	       
\newtheorem{example}[thm]{Example}

\makeatletter
\def\theequation{\@arabic\c@equation}
\def\thethm{\@arabic\c@thm}
\def\thelem{\@arabic\c@thm}
\def\thecrlr{\@arabic\c@thm}
\def\theprp{\@arabic\c@thm}
\def\therem{\@arabic\c@thm}
\def\theexample{\@arabic\c@thm}
\makeatother

\begin{document}


\title[Non-linear systems in divergence form] {Boundedness of the  solutions to  nonlinear  systems with  Morrey data}

\author[L.G. Softova]{Lubomira G. Softova}

\address{Department of Civil Engineering,\\
Design, Construction and Environment\\
  Second University  of Naples\\
 Via Roma 29\\
 81031 Aversa\\
 Italy}

\email{luba.softova@@unina2.it}



\maketitle

\begin{abstract}
We consider nonlinear elliptic systems satisfying componentwise coercivity condition. The nonlinear terms have controlled growths with respect to the solution and its gradient, while the behaviour in the independent variable is governed by functions in Morrey spaces. We firstly prove essential boundedness of the weak solution and then obtain Morrey regularity of its gradient.  
\end{abstract}

\section{Introduction}

Let $\Omega\subset\R^n, n\geq 2$ be a bounded domain satisfying the \eqref{A}-condition. 
We are interested in  boundedness and Morrey regularity of the weak solutions to nonlinear
 elliptic systems of the type
\begin{equation}\label{NS0}
-\div\bA(x,\bu, D\bu)+\bb(x,\bu,D\bu)=\ff(x),\qquad x\in \Omega
\end{equation}
where the nonlinear terms are Carath\'eodory maps 
\begin{align*}
\bA(x,\bu,\bz): &\ \Omega\times\R^N\times\M^{N\times n}\to \R^{N\times n},\\
\bb(x,\bu,\bz):  &\  \Omega\times\R^N\times\M^{N\times n}\to \R^{N}.
\end{align*}

The celebrated   result of    De Giorgi  \cite{DeG} and  Nash \cite{N} implies   that any weak solution  $u\in  W_0^{1,2}(\Omega)$  of the linear elliptic equation   
 $D_i(A_{ij}(x)D_ju +g_i(x))=f(x)$ 
is locally H\"older continuous when $g_i\in L^p$ with $p>n$ and $f\in L^q$ with $q>n/2,$ even if the coefficients are only $L^\infty.$  Unfortunately 
the  De Giorgi-Nash result  does not hold anymore  if we consider a system of uniformly elliptic equations  because of the lack of {\it Maximum principle}. This  was shown by De Giorgi himself almost ten years later, constructing  a counterexample  \cite{DG}. 
Precisely, the function  $\bu=1-x/|x|^\gamma\in W_0^{1,2}(\B_1(0);\R^n)$ is a solution to
$$
D_i(A_{ij}^{\alpha\beta}(x)D_ju^\beta(x))=0\qquad \text{ in }  \B_1(0)
$$
with suitably chosen coefficients $A_{ij}^{\alpha\beta}\in L^\infty(\B_1(0)).$

Moreover, the result of De Giorgi-Nash cannot be extended to quasilinear systems even if the coefficients are analytic functions, as it was  shown by Giusti  and  Miranda in  \cite{GMir}.
In order to get a maximum principle for elliptic systems we need  to impose some quite restrictive structural conditions. The simplest one requires the   system to be  in diagonal form, or  {\it decoupled.}
\begin{example}
Consider the  operator 
$
\div (\bA(x,D\bu))=0$ in $\Omega
$
with coefficients  
$$
A_i^\alpha(x,D\bu)=\sum_{j=1}^n\sum_{\beta=1}^N\delta_{\alpha\beta}A^{\alpha\beta}_{ij}(x)D_ju^\beta
$$  
where $\delta_{\alpha\beta}$ is the Kronecker delta.  
Then  $u^\alpha$ solves a single elliptic equation and $
\sup_\Omega u^\alpha\leq \sup_{\partial \Omega} u^\alpha,
$   for each  $\alpha=1,\ldots,N.$ 
\end{example}
One more example was given by Ne\v{c}as  and  Star\'a in \cite{NS}.
\begin{example}  Consider the system
$ \div\bA(x,\bu,D\bu)=0$   in $ \Omega $
 that is diagonal for   large values of $u^\alpha,$ that is,
\begin{equation}\label{exNS}
0<\theta_\alpha\leq u^\alpha\  \Longrightarrow\   A_i^\alpha(x,\bu,D\bu)= \sum_{j=1}^n\sum_{\beta=1}^N\delta_{\alpha\beta} 
A_{ij}^{\alpha\beta}(x,\bu)D_ju^\beta
\end{equation}
with bounded and elliptic $A_{ij}^{\alpha\beta}.$ It turns out that
$$
\sup_\Omega u^\alpha\leq \max\big\{ \theta_\alpha; \sup_{\partial\Omega} u^\alpha  \big\}
$$
also in this case.
\end{example}

The situation becomes more complicated if we consider {\it general nonlinear  system}
\begin{equation}\label{NS}
\div\bA(x,\bu, D\bu)=\bb(x,\bu,D\bu)\,.
\end{equation}
Along with the Carath\'eodory conditions on the maps $\bA(x,\bu,\bz)$ and $\bb(x,\bu,\bz)$ we need to control also the growths of $\bA$ and $\bb$ with respect to $\bu$ and $\bz.$ These additional {\it controlled growth conditions} ensure the convergence of the integrals in the definition of {\it weak solution} to \eqref{NS} (see \eqref{weak}).

In \cite{LP} Leonetti and Petricca  assume  {\it componentwise  coercivity condition}  on $\bA$ and positivity of $\bb$ for large values of $u^\alpha,$ that is, there exist positive constants $\theta_\alpha$ such that
\begin{equation}\label{exLP}
\theta_\alpha\leq u^\alpha\quad  \Longrightarrow\quad 
\begin{cases}
\ds   \nu|\bz^\alpha|^p -M_\alpha\leq \sum_{i=1}^n A_i^\alpha(x,\bu,\bz)z_i^\alpha\\
 \ds  0\leq b^\alpha(x,\bu,\bz)\,.
\end{cases}
\end{equation}
Combining the {\it Sobolev   inequality} with the {\it Stampacchia Lemma \cite{Stamp}} they  get a componentwise bound of the solution, covering this way also the systems studied in \cite{NS},  since \eqref{exNS} is a special case of \eqref{exLP}. 
Let us note that getting essential boundedness of the weak solution  to \eqref{NS0} is a starting point for a further  study of its regularity  in various function spaces. In \cite{DK, P,  PS2}  the authors  obtain better integrability and H\"older regularity of the bounded solutions to quasilinear elliptic equations $(N=1)$ under controlled growth conditions on the nonlinear terms. Further this result has been extended in \cite{Sf} to semilinear  uniformly elliptic  systems of the form 
\begin{equation}\label{QLS}
\div(\bA(x)D\bu+\ba(x,\bu))=\bb(x,\bu,D\bu)\qquad \text{ in }  \Omega
\end{equation}
with minimal regular assumptions on the coefficients and the underlying domain.  Precisely,  it is  shown  that if  the nonlinear terms 
satisfy the controlled growth conditions \eqref{contr1} with $\varphi\in L^p(\Omega),$ $p>2$ and $\psi\in L^q(\Omega), $ 
$q>\frac{2n}{n+2}$ then any bounded weak solution to \eqref{QLS}  belongs to $W_0^{1,r}$ with $r=\min\{p,q^*\}.$ 

The natural question that arises is what kind of regularity of the solution to \eqref{NS0}  we can  expect if the given functions $\varphi$ and $\psi$ belong to some   Morrey spaces. In the case of a single equation we count  with the result of Byun and Palagachev \cite{BP}.
 Combining the Gehring-Giaquinta-Modica  lemma, the Adams trace inequality and the Hartmann-Stampacchia maximum principle they obtain $L^\infty$ estimate of  the solution.  Further, the  Morrey-type estimate of the gradient  permits the authors  to show also H\"older regularity of the solution.

Our goal is to obtain a  componentwise maximum principle for any  solution of \eqref{NS} supposing that the operators $\bA$ and $\bb$ satisfy structural conditions expressed in terms of Morrey functions. As a consequence we obtain also Morrey regularity of the gradient of $\bu$ extending  such a way  the regularity  results obtained in  \cite{BP, DK,  LP, NS,   PSesyst, Sf} to nonlinear systems with Morrey data.

 Recall that a  real valued function $f\in L^p(\Omega)$  belongs  to the Morrey space $L^{p,\lambda}(\Omega)$ with $p\in[1,\infty),$ $\lambda\in(0,n),$   if  
\begin{equation}\label{defMorrey}
 \|f\|_{p,\lambda;\Omega}=\left( \sup_{\B_r(x) }\frac1{r^\lambda}
\int_{\Omega\cap\B_r(x)}|f(y)|^p\,dy   \right)^{1/p}<\infty
\end{equation}
where the supremum is taken over all balls $\B_r(x),$   $r\in(0,\diam\Omega]$ and $x\in \ol\Omega.$    Working in the framework of the Morrey spaces we note that the Sobolev trace inequality is not enough anymore. For this goal we will use the following result due to Adams. 
\begin{lem}[Adams Trace Inequality,  \cite{AdTrace, ChTrace, Rok}]\label{AdamsTrace}
Let $m$ be a positive Radon measure with support in $\Omega$ and such that for each ball $\B_\rho$ it holds 
\begin{equation}\label{tau0}
m(\B_{\rho}) \leq K \rho^{\tau_0},\quad 
\tau_0=\frac{s}{r}(n-r),\quad 1<r<s<\infty, \quad r<n
\end{equation}
with an absolute constant $K>0.$
Then 
\begin{equation}\label{Adams1}
\left(\int_\Omega |v(x)|^s\, dm \right)^{\frac1s}\leq C(n,s,r)K^{\frac1s}\left( \int_\Omega |Dv(x)|^r\, dx \right)^{\frac1r} 
\end{equation}
for each  function  $ v\in W_0^{1,r}(\Omega).$
\end{lem}

In what follows we suppose that $\Omega\subset \R^n, n\geq2,$ is a bounded domain satisfying the \eqref{A}-condition,  that is, 
 there exists   a  constant $A_\Omega>0$  such that 
\begin{equation}\tag{A}\label{A}
|\Omega_r(x)|\geq A_\Omega\,  r^n\qquad \forall\  x\in\ol\Omega, \  r\in(0,\diam\Omega]
\end{equation}
where 
$\Omega_r(x)=\Omega\cap \B_r(x).$ It is worth noting that the \eqref{A}-condition 
excludes interior cusps at each point of the boundary and  guarantees
the validity of the Sobolev embedding theorem in $W^{1,p}(\Omega).$ This geometric property is surely satisfied when $\partial\Omega$ has the uniform interior cone property (e.g. $C^1$-smooth or Lipschitz continuous boundaries), but it holds also for the
 Reifenberg falt domains 
boundaries (cf. \cite{PS2}).

   Throughout  the text  the standard summation convention on the repeated  indexes is adopted. 
The letter $C>0$ is used for various constants and may change from one occurrence to another.

\section{Maximum principle}

Consider the nonlinear system
\begin{equation}\label{NonlSyst}
- D_i\big(A_{i}^{\alpha}(x,\bu, D\bu)\big)+ b^\alpha(x,\bu,D\bu) =f^\alpha(x) \  \text{ in } \Omega
\end{equation}
where $
\bA=\{A^\alpha_i(x,\bu,\bz)\}_{i\leq n}^{\alpha\leq N}
$ and 
$  
\bb=(b^1(x,\bu,\bz),\ldots, b^N(x,\bu,\bz))
$
 are  measurable in $x\in\Omega$ and continuous in $(\bu,\bz)$ for almost all (a.a.) $x\in\Omega.$ Suppose that for each 
 $(x,\bu,\bz)\in\Omega\times \R^N\times  \Mb^{N\times n}$  the following   {\it controlled growth conditions} hold. Namely,  
\begin{equation}\label{contr1}
\begin{cases}
\ds  |\bA(x,\bu,\bz) |\leq \Lambda\big(\varphi(x)+|\bu|^{\frac{\nu}2}+|\bz|\big)\\[5pt]
\ds   |\bb (x,\bu,\bz)|\leq \Lambda\big(\psi(x) + |\bu|^{\nu-1}+|\bz|^{2\frac{\nu-1}{\nu})}  \big) 
\end{cases}
\end{equation}
as   $  |\bu|, |\bz|\to\infty,$
with some positive constant $\Lambda.$  Here
$\nu$ is the Sobolev conjugate of $2,$ that is, 
\begin{equation}\label{conjugate}
\nu=\begin{cases}
\ds \frac{2n}{n-2} & \text{  if } n\geq 3\\
\text{any   number } > 2 & \text{  if  } n =2,
\end{cases}
\end{equation}
and the given functions $\varphi,$  $\psi$ and $f^\alpha$ satisfy
\begin{equation}\label{regularity}
\begin{cases}
\ds \varphi\in L^{p,\lambda}(\Omega), & p>2, \   \lambda\in(0,n), \   p+\lambda>n\\[5pt]
\ds \psi,  f^\alpha \in L^{q,\mu}(\Omega),   & q>\frac{\nu}{\nu-1}, \  \mu\in(0,n), \  2q+\mu>n.
\end{cases}
\end{equation}

In the particular case  $n=2$ the powers of $|\bu|$ could be arbitrary positive numbers  greater then 1,  while the growth of $|\bz|$ is strictly sub-quadratic (cf. \cite{G, LU}).

Under  a {\it weak solution} of  \eqref{NonlSyst} we mean a function  
 $\bu\in W^{1,2}(\Omega;\R^N)\cap L^\nu(\Omega;\R^N),$   satisfying
\begin{align}\label{weak}
\nonumber
\sum_{i=1}^n  \int_\Omega  A^\alpha_i(x,\bu, D\bu)D_i\phi^\alpha(x)\, dx &  
 +\int_\Omega b^\alpha(x,\bu(x),D \bu(x))\phi^\alpha(x)\, dx\\
& = \int_{\Omega} f^\alpha(x) \phi^\alpha(x)\,dx
\end{align}
for all $\phi=(\phi^1,\ldots,\phi^N)\in W_0^{1,2}(\Omega;\R^N).$   
The conditions  \eqref{contr1}-\eqref{regularity}
are the natural ones that ensure the convergence of the integrals in \eqref{weak}. Moreover, they are optimal  as it is seen from the following example in the case of  single equation (cf. \cite{LSU, P}). 
\begin{example}
The function $
u(x)=|x|^{\frac{r-2}{r-1}}\in W^{1,2}(\B_1(0)),$ with  $ n\geq3 $ and  $  \frac{n+2}{n}<r<2$
  is a solution  to the equation $\Delta u=C|Du|^r$  in $\B_1(0).$ Note that $u\not\in L^\infty(\B_1(0)).$ 
\end{example}

Generally we cannot expect boundedness of the  solutions  to \eqref{NonlSyst} unless we add some restrictions on the  structure of the operator (see for example \cite{JS, LP}). 
For this goal  we impose componentwise  coercivity on  $A_i^\alpha$  and a  sign condition on $b^\alpha.$ 

For  every $\alpha\in\{1,\ldots,N\}$ there exist positive constants  $\theta_\alpha,$ $ \gamma $ and a function $ \varphi $ such that  for each  $u^\alpha\geq \theta_\alpha$  we have
\begin{equation}\label{coercivity}
\begin{cases}
\ds \gamma|\bz^\alpha|^2 - \Lambda \varphi(x)^2\leq \sum_{i=1}^n A_i^\alpha (x,\bu,\bz)z_i^\alpha\\
\ds\varphi\in L^{p,\lambda}(\Omega), \ p>2, \ p+\lambda>n\\
  \ds 0\leq b^\alpha(x,\bu,\bz)\quad  \text{ for a.a. } x\in\Omega, \  \forall  \  \bz\in \Mb^{N\times n}\,.
\end{cases}
\end{equation}

\begin{thm}[Maximum principle]\label{thm1}
Let $\Omega$ be \eqref{A}-type  domain  and  $\bu\in W^{1,2}(\Omega;\R^N)\cap L^\nu(\Omega; \R^N)$ be a weak solution to  \eqref{NonlSyst} under the conditions   \eqref{contr1}, \eqref{regularity}  and  \eqref{coercivity} and such that $\sup_{\partial\Omega}u^\alpha<\infty$.  Then  
$$
\sup_\Omega u^\alpha \leq\max\{ \theta_\alpha, \sup_{\partial\Omega} u^\alpha \}+ M_\alpha\qquad  \alpha\in \{1,\ldots,N\}
$$
where $M_\alpha$ depends on $ n, p, \lambda,  \Lambda, \gamma,$ $\|\varphi \|_{p,\lambda;\Omega},$  $\| f^\alpha \|_{q,\mu;\Omega}$  and $  |\Omega|.$ 
\end{thm}
\begin{pf}
We choose  a constant $L>0$ such that 
$
L\geq \max\{ \theta_\alpha; \sup_{\partial\Omega}u^\alpha \}
$
and define the set
$
\cA_L^\alpha=\{ x\in\Omega:\  u^\alpha(x)-L>0  \}\,.
$
Then we take  a  vector  function $\bv$
 as follows 
$$
v^\beta=
\begin{cases}
 \ds \max\,\{u^\alpha-L; 0\} &\text{if } \beta=\alpha\\
0 &\text{if } \beta\not=\alpha
\end{cases},\qquad D v^{\beta}=\begin{cases}
  D u^\alpha\chi_{\cA_L^\alpha} &\text{if } \beta=\alpha\\
0 &\text{if } \beta\not=\alpha
\end{cases}.
$$
It is clear that $
\bv \in W_0^{1,2}(\Omega;\R^N)$ and hence $\bv\in L^\nu(\Omega;\R^N)$ by the Sobolev embedding.
Choosing $\phi^\alpha=v^\alpha$ as a test function we obtain 
\begin{align*}
 \int_{\cA^\alpha_L} A_i^\alpha(x,\bu,D\bu)D_iu^\alpha(x) \,dx \  +& \int_{\cA^\alpha_L} b^\alpha(x,\bu,D\bu)(u^\alpha(x)-L)\,dx\\
=& \int_{\cA^\alpha_L}f^\alpha(x)(u^\alpha(x)-L)\,dx\,.
\end{align*}

We start with  the case
 $ n\geq 3$ when $\nu=2n/(n-2).$ 
Define  the Radon measure $dm$ supported in $\Omega$ by
$$
dm:=(\chi_\Omega(x) +\varphi(x)^2+|f^\alpha(x)|)\,dx,
$$
where $\chi_\Omega$ is the characteristic function of $\Omega.$
Then by   \eqref{coercivity} we get the estimate 
\begin{align}\label{eq1}
\nonumber
 \int_{\cA^\alpha_L}   |Du^\alpha(x)|^2\,dx\leq &\    \frac{\Lambda}{\gamma}\int_{\cA^\alpha_L} \varphi(x)^2\,dx+\frac1{\gamma}
\int_{\cA^\alpha_L}|f^\alpha(x)|(u^\alpha(x)-L)\,dx\\
\nonumber
\leq   &\     \frac{\Lambda}{\gamma} \int_{\cA^\alpha_L} (\chi_\Omega(x)
  +\varphi(x)^2 +|f^\alpha(x)|)\,dx\\
&\  +  \frac{1}{\gamma} \int_{\cA^\alpha_L}|f^\alpha(x)|(u^\alpha(x)-L)\,dx\\
\nonumber
\leq &\   C(\Lambda,\gamma)\big( m(\cA^\alpha_L )+  J \big).
\end{align}
In order to estimate the integral $J= \int_{\cA^\alpha_L}|f^\alpha(x)|(u^\alpha(x)-L)\,dx $ we  make use of the Lemma~\ref{AdamsTrace} applied to  the Radon measure $dm'=|f^\alpha(x)|dx.$ Hence
$$ 
J=\int_{\cA^\alpha_L}(u^\alpha(x)-L)dm'\leq \left(\int_{\cA_L^\alpha}|u^\alpha(x)-L|^{s'}\,dm'  \right)^{\frac1{s'}} m'(\cA^\alpha_L)^{1-\frac1{s'}}\,.
$$
Evaluating the measure $m'$ over  a ball $\B_\rho$ we get
\begin{align*}
m'(\B_\rho)=& \  \int_{\B_\rho} |f^\alpha(x)|\,dx  \leq C(n)\rho^{n-\frac{n-\mu}q}\left(\frac1{\rho^\mu}
 \int_{\B_\rho} |f^\alpha(x)|^q\,dx  \right)^{1/q}\\
\leq&\  C(n)\rho^{n-\frac{n-\mu}q} \|f^\alpha\|_{q,\mu;\Omega}=K \rho^{n-\frac{n-\mu}q}
\end{align*}
with $K=K(n,q,\diam\Omega,   \|f^\alpha\|_{q,\mu;\Omega} ).$
We apply now the Lemma~\ref{AdamsTrace} with $r'=2,$ $\tau'_0=n-\frac{n-\mu}{q}$ and
 $s'=\frac2{n-2}\big(n-\frac{n-\mu}q  \big)>2,$ calculated via \eqref{tau0}. Hence 
\begin{align}\label{eq2}
\nonumber
J\leq &\  C K^{\frac1{s'}}\left(  \int_{\cA^\alpha_L} |Du^\alpha(x)|^2\,dx\right)^{\frac12} m'(\cA^\alpha_L)^{1-{\frac1{s'}}}\\
\leq &\   C\left[ \varepsilon \int_{\cA_L^\alpha} |Du^\alpha(x)|^2\,dx+\frac1{\varepsilon}
 m'(\cA_L^\alpha)^{2(1- {\frac1{s'}})}  \right]\,.
\end{align}
Combining \eqref{eq1} and \eqref{eq2}, taking $\varepsilon $  small enough,  moving the integral of the gradient on the left-hand side, and keeping in mind that $2(1-\frac1{s'})>1$ and $m'(\cA^\alpha_L)\leq m(\cA^\alpha_L)$  we obtain
\begin{equation}\label{eq3}
\int_{\cA^\alpha_L}|Du^\alpha(x)|^2\,dx\leq C\big( m(\cA^\alpha_L) +m(\cA^\alpha_L)^{2(1-\frac1{s'})} \big)\leq
C m(\cA^\alpha_L)
\end{equation}
where  the constant $C$ depends on known quantities.

To complete the estimate \eqref{eq1} we will use once again the  Lemma~\ref{AdamsTrace}. It is immediate that
$m(\B_\rho)$ of a ball
 $\B_\rho\Subset\Omega$ is
\begin{align}\label{measure_bound}
\nonumber
m(\B_{\rho}) & =   \int_{\B_\rho}\big(\chi_\Omega(x)+\varphi(x)^2+|f^\alpha(x)|   \big)\,dx\\
&\leq C(n) \rho^n+ \rho^{n-\frac{2(n-\lambda)}p} \|\varphi\|^2_{p,\lambda;\Omega}+ 
\rho^{n-\frac{n-\mu}{q}}\|f^\alpha\|_{q,\mu;\Omega}
\leq K \rho^{\tau_0}
\end{align}
with $K=K(n,p,q,\diam \Omega, \|\varphi\|_{p,\lambda;\Omega},\|f^\alpha\|_{q,\mu;\Omega} )$ and 
$$
\tau_0=\min\left\{n-\frac{2(n-\lambda)}p; n-\frac{n-\mu}q    \right\} >n-2.
$$
Applying \eqref{Adams1} with $r=2<n$ and calculating  $s$ from  \eqref{tau0} 
  we get
\begin{align}\label{estimate1}
\nonumber
\int_{\cA^\alpha_L} (u^\alpha(x)-L)\,dm &\leq  \left( \int_{\cA^\alpha_L}|u^\alpha(x)-L|^s\,dm  \right)^{\frac1s} m(\cA^\alpha_L)^{1-\frac1s}\\
&\leq  C K^{\frac1s}
\left(  \int_{\cA^\alpha_L}|D u^\alpha(x)|^2\,dx \right)^{\frac12}   m(\cA^\alpha_L)^{1-\frac1s}   \\
\nonumber
&\leq C(n,p,K,\gamma,\Lambda)m(\cA^\alpha_L)^{1+\frac12-\frac1s}
\end{align}
with $s=\min\left\{\frac{2np-4(n-\lambda)}{p(n-2)} ; \frac{2nq-2(n-\mu)}{q(n-2)}  \right\}>2.$ 

\smallskip

A similar bound holds also in the case $n=2.$ In fact, for any ball 
 $\B_\rho\subset\R^2$ we have
$$
m(\B_\rho)\leq C\rho^2+  \rho^{2-\frac{2(2-\lambda)}p}  \|\varphi\|^2_{p,\lambda;\Omega} +
 \rho^{2-\frac{2-\mu}{q}}  \|f^\alpha\|_{q,\mu;\Omega} \leq K \rho^{\tau_0}
$$
with $\tau_0=\min\left\{ 2-\frac{2(2-\lambda)}{p}; 2-\frac{2-\mu}{q}  \right\}>0.$
Choosing $s=2$ we calculate $r$ from  \eqref{tau0}  
$$
r=\max\left\{  \frac{2p}{2p-2+\lambda}; \frac{4q}{4q-2+\mu} \right\}\in(1,2).
$$
Then by the H\"older and the Adams trace inequalities we obtain 
\begin{align}\label{estimate2}
\nonumber
\int_{\cA^\alpha_L} (u^\alpha(x)&-L)\,dm   \leq\left(\int_{\cA^\alpha_L} (u^\alpha(x)-L)^2\,dm \right)^{\frac12}m(\cA^\alpha_L)^{\frac12}\\
& \leq C K^{\frac12}\left( \int_{\cA^\alpha_L} |Du^\alpha(x)|^r\,dx \right)^{\frac1r} m(\cA^\alpha_L)^{\frac12}\\
\nonumber
& \leq C K^{\frac12}\left( \int_{\cA^\alpha_L}|Du^\alpha(x)|^{2}\,dx \right)^{\frac{1}2} 
 \left(\int_{\cA^\alpha_L}\chi_\Omega(x)\, dx \right)^{\frac1r-\frac{1}2}m(\cA^\alpha_L)^{\frac12}\\
\nonumber
& = CK^{\frac12}  \left( \int_{\cA^\alpha_L}|Du^\alpha(x)|^{2}\,dx \right)^{\frac{1}2}   m(\cA^\alpha_L)^{\frac1r}\,.
\end{align}
In order to estimate the integral in the last term  we go back to \eqref{eq1}. Consider again the Radon measure 
$dm'=|f^\alpha(x)|dx$ and calculate $m'(\B_\rho)\leq  K\rho^{2-\frac{2-\mu}{q}}.$  Then choosing $s'=2$ we  get $r'=\frac{4q}{4q-(2-\mu)}\in(1,2)$ from \eqref{tau0}. This way, the  Lemma~\ref{AdamsTrace} and the H\"older inequality give
\begin{align}\label{eq4}
\nonumber
J\leq &\  \left(\int_{\cA_L^\alpha}|u^\alpha(x)-L|^2\,dm'  \right)^{\frac12} m'(\cA^\alpha_L)^{\frac12}\\
\leq &\   C K^{\frac12}  \left(\int_{\cA_L^\alpha}|Du^\alpha(x)|^{r'}\,dx  \right)^{\frac1{r'}} m'(\cA^\alpha_L)^{\frac12}\\
\nonumber
\leq &\   C K^{\frac12}   \left(\int_{\cA_L^\alpha}|Du^\alpha(x)|^2\,dx  \right)^{\frac12}   m(\cA^\alpha_L)^{\frac1{r'}-\frac12}  m(\cA^\alpha_L)^{\frac12}\\
\nonumber
\leq &\  C\left[\varepsilon \int_{\cA_L^\alpha}|Du^\alpha(x)|^2\,dx +\frac1{\varepsilon} m(\cA_L^\alpha)^{\frac2{r'}}     \right]\,.
\end{align}
Unifying \eqref{eq1} and \eqref{eq4},   taking $\varepsilon$ small enough and keeping in mind that $\frac2{r'}>1,$ we get
$$
\int_{\cA_L^\alpha}|Du^\alpha(x)|^2\,dx\leq C m(\cA_L^\alpha)
$$ 
where the constant depends on the same  quantities as in \eqref{eq3}. Then the estimate \eqref{estimate2} becomes 
\begin{equation}\label{eq5}
\int_{\cA^\alpha_L}(u^\alpha(x)-L)\,dm\leq C m(\cA^\alpha_L)^{1+\frac1r-\frac12}.
\end{equation}
Unifying the estimates \eqref{estimate1} and \eqref{eq5} we obtain
\begin{equation}\label{estimate3}
\int_{\cA^\alpha_L}(u^\alpha(x)-L)\, dm\leq C m(\cA^\alpha_L)^{1+\sigma_0}
\end{equation}
where 
$$
\sigma_0=\begin{cases}
\ds\frac12-\frac1s=\max\big\{\frac{p+\lambda-n}{np-2(n-\lambda)};\frac{\mu+2q-n}{2nq-2(n-\mu)}   \big\} & \text{ if } n>2\\[10pt]
 \ds\frac1r-\frac12= \min\big\{\frac{p+\lambda-2}{2p};  \frac{2q+\mu-2}{4q} \big\}  & \text{ if } n=2\,.
\end{cases}
$$

Suppose now that 
$
m(\cA^\alpha_L)>0,
$
otherwise $\sup_{\Omega}u^\alpha(x)\leq L.$    For any $L_1>L$ we have  $\cA^\alpha_{L_1}\subset \cA^\alpha_L$ and therefore \eqref{estimate3} yields
\begin{align*}
(L_1-L)m(\cA^\alpha_{L_1})&\leq \int_{\cA^\alpha_{L_1}}(u^\alpha(x)-L)\,dm\\
&\leq \int_{\cA^\alpha_L}(u^\alpha(x)-L)\,dm\leq C m(\cA^\alpha_L)^{1+\sigma_0}\,.
\end{align*}
Hence 
$$
m(\cA^\alpha_{L_1})\leq \frac{C}{L_1-L}m(\cA^\alpha_L)^{1+\sigma_0}\,.
$$
In order to estimate the measure of the set $\cA_L^\alpha$ we will apply  the following 
Maximum Principle due to
 Stampacchia  \cite[Lemma~4.1]{Stamp}.
\begin{lem}\label{lemMaximum}
Let $\Theta:[L_0,\infty)\to [0,\infty)$ be a decreasing function. Assume that there exist $c,a\in (0,\infty)$ and $b\in (1,\infty)$
such that 
$$
L_1>L\geq L_0 \  \Longrightarrow \  \Theta(L_1)\leq \frac{c}{(L_1-L)^\alpha}(\Theta(L))^b.
$$
Then 
$$
\Theta(L_0+d)=0 \quad \text{ where }\quad d=\big[c\Theta(L_0)^{b-1} 2^{\frac{ab}{b-1}}  \big]^{\frac1a}\,.
$$
\end{lem}

The application of the  Lemma~\ref{lemMaximum} to the function  $\Theta(L)=m(\cA^\alpha_L)$ with    $a=1,$  $b=1+\sigma_0$ and   $L_0=\max\{\theta_\alpha,\sup_{\partial\Omega} u^\alpha \}$  yields
\begin{equation}\label{estimate4}
m(\cA^\alpha_{L_0+d_\alpha})=0 \qquad \text{ where } \qquad d_\alpha\leq C m(\Omega)^{\sigma_0} 2^{1+\frac1{\sigma_0}}\,.
\end{equation}

The last assertion means that for each $\alpha=1,\ldots,N$ there exists a constant $M_\alpha$ depending on $n,$ $p,$ $\lambda,$ $q,$ $\mu,$  $\gamma,$  $\Lambda,$ $|\Omega|,$  $\|\varphi\|_{p,\lambda;\Omega}$ and 
$\|f^\alpha\|_{q,\mu;\Omega}$ such that 
\begin{equation}\label{Maximum}
\sup_{\Omega}u^\alpha<\max\{\theta_\alpha; \sup_{\partial\Omega}u^\alpha\}+M_\alpha
\end{equation}
and this completes the proof of Theorem \ref{thm1}
\end{pf}

\section{The Dirichlet Problem}

We study the boundedness and the Morrey regularity  of the weak solutions to the following Dirichlet problem
\begin{equation}\label{DP3}
\begin{cases}
- \div\bA(x,\bu(x) D\bu(x))+ \bb(x,\bu,D\bu) =\ff(x) & x\in\Omega \\
\bu(x)=0\quad  & x\in  \partial\Omega
\end{cases}
\end{equation}
in a bounded domain $\Omega\subset\R^n.$
\begin{thm}[Essential Boundedness of the Solution]\label{crlr1}
Let $\bu\in W_0^{1,2}(\Omega,\R^N)$ be a solution to \eqref{DP3}  and assume  \eqref{A}, \eqref{contr1}, and \eqref{regularity}. Suppose in addition that
\begin{equation}\label{coercivity2}
\begin{cases}
 \ds  \gamma |\bz^\alpha|^2-\Lambda \varphi(x)^2\leq \sum_{i=1}^n A_i^\alpha(x,\bu,\bz)z_i^\alpha\\[6pt]
\ds  \varphi(x)\in L^{p,\lambda}(\Omega), \  p>2, \  \lambda\in(0,n), \   \lambda+p>n,\\[6pt]
\ds 0\leq   b^\alpha(x,\bu,\bz)\, \sign  u^\alpha(x)
\end{cases}
\end{equation}
for $|u^\alpha|\geq \theta_\alpha>0.$ 
Then there exists a constant $M$ depending on known quantities  such that 
$$
\|\bu\|_{\infty;\Omega}\leq M\,.
$$
\end{thm}
\begin{pf}
Take a positive constant $L$ such that $L\geq \theta_\alpha$ and
consider  the  set $\bar \cA^\alpha_L=\{x\in\Omega: u^\alpha(x)+L<0  \}.$  Then the  Theorem~\ref{thm1} applied  to $-u^\alpha$ gives
\begin{equation}\label{eq6}
\inf_\Omega u^\alpha>-\theta_\alpha - M_\alpha\,.
\end{equation}
Unifying \eqref{Maximum} and \eqref{eq6} we get boundedness of $\|u^\alpha\|_{\infty;\Omega}$ for each $\alpha=1,\ldots,N.$  Then  
$$
\|\bu\|_{\infty;\Omega}=\max_{1\leq \alpha\leq N}\|u^\alpha\|_{\infty;\Omega}=: M<\infty\,.
$$

\end{pf}

\begin{thm}[Morrey regularity of the gradient]\label{MorreyEst}
Let $\Omega $ be a bounded \eqref{A}-type domain in $\R^n, n\geq3,$ and  $\bu\in  W_0^{1,2}(\Omega,\R^N) $ be a weak solution to \eqref{DP3} under the assumptions  \eqref{contr1}, \eqref{regularity}, and \eqref{coercivity2}. Then  $D\bu\in L^{2,n-2}(\Omega,\R^N)$ and
\begin{equation}\label{Gradient4}
\int_{\Omega_\rho(x)}|Du(y)|^2\,dy\leq C \rho^{n-2}\qquad \forall x\in\Omega, \  \rho\in(0,\diam\Omega]
\end{equation}
with a constant depending on known quantities.
\end{thm}
\begin{pf}
Fix $x_0\in\Omega$ and $\rho>0$ be such that $\B_\rho(x_0)\subset \B_{2\rho}(x_0)\Subset \Omega,$  $\rho>1.$ Define a cut-off function $\zeta(x)\in C^1(\R^n)$  
$$
\zeta(x)=\begin{cases}
1  & x\in \B_{\rho}(x_0),\\
0 & x\not\in \B_{2\rho}(x_0),
\end{cases}\qquad |D\zeta|\leq \frac{C}{\rho}\,.
$$ 
For any fixed $\alpha$  take  $\phi^\alpha(x)=e^{u^\alpha(x)}\zeta(x)^2  $ as a test function in \eqref{weak} to get
\begin{align*}
  \sum_{i=1}^n \int_{\Omega}&   A_i^\alpha(x,\bu, D\bu)   e^{u^\alpha(x)}D_iu^\alpha(x) \zeta(x)^2\,dx\\ 
&= \int_{\Omega}  f^\alpha(x) e^{u^\alpha(x)}2\zeta(x) D_i\zeta(x)\,dx\\
&\quad -  \sum_{i=1}^n \int_{\Omega} A_i^\alpha(x,\bu, D\bu)e^{u^\alpha(x)}2\zeta(x) D_i\zeta(x)\,dx\\
&\quad -\int_{\Omega}b^\alpha(x,\bu,D\bu)e^{u^\alpha(x)}\zeta(x)^2\,dx\,.
\end{align*}
The left-hand side can be estimated by \eqref{coercivity2} while for the right-hand side we use   \eqref{contr1} and \eqref{regularity}
\begin{align*}
 \sum_{i=1}^n& e^{-M}\int_{\Omega}
\big(\gamma |D u^\alpha(x)|^2 -\Lambda \varphi(x)^2  \big)\zeta(x)^2 \,dx\\
 &  \leq 2e^{M}  \int_{\Omega} |  f^\alpha(x)| \zeta(x) |D\zeta(x)|\,dx\\
&\quad + 2 n \Lambda  e^{M}\int_{\Omega} \big(\varphi(x)+|\bu|^{\frac{n}{n-2}} +|D\bu|  \big)\zeta(x) |D\zeta|\,dx\\
 &\quad +  \Lambda e^{M} \int_{\Omega}\big( \psi(x)+|\bu|^{\frac{n+2}{n-2} }+|D\bu|^{\frac{n+2}{n}} \big)  \zeta(x)^2\,dx.
\end{align*}
To proceed further, we use the Young inequality
$
ab\leq  \varepsilon a^p+\frac{b^{p/(p-1)}}{\varepsilon^{1/(p-1)}},
$   whence
\begin{align*}
\int_{\Omega}|f^\alpha(x)| \zeta(x) |D\zeta(x)|\,dx\leq &\ \frac12\int_{\Omega}| f^\alpha(x)|^2 \zeta(x)^2\,dx 
+\frac12\int_{\Omega}|D\zeta(x)|^2\,dx\\
\int_{\Omega}\varphi(x) \zeta(x) |D\zeta(x)|\,dx\leq&\ \frac12\int_{\Omega} \varphi(x)^2 \zeta(x)^2\,dx 
+\frac12\int_{\Omega}|D\zeta(x)|^2\,dx\\
\int_{\Omega} |\bu|^{\frac{n}{n-2}}\zeta(x) |D\zeta(x)|\,dx \leq&\ \frac12 M^{\frac{n}{n-2}}\left( \int_{\Omega}\zeta(x)^2\,dx+ \int_{\Omega}|D\zeta(x)|^2\,dx \right)\\
\int_{\Omega}|D\bu|\zeta(x) |D\zeta|\,dx\leq&\ \varepsilon\int_{\Omega}|D\bu|^2\zeta(x)^2\,dx +\frac1{\varepsilon} \int_{\Omega}|D\zeta(x)|^2\,dx\\
\int_{\Omega} |D\bu|^{\frac{n+2}{n}}\zeta(x)^2\,dx\leq&\ \varepsilon\int_{\Omega} |D\bu|^2\zeta(x)^2\,dx
+\varepsilon^{-\frac{n+2}{n-2}}\int_{\Omega}\zeta(x)^2\,dx\,.
\end{align*}
Unifying the above estimates we get 
\begin{align}\label{Gradient}
\nonumber
\int_{\Omega}& |Du^\alpha(x)|^2\zeta(x)^2\,dx\\
 &\ \leq C\int_{\Omega}(1+|f^\alpha(x)|+\psi(x)+\varphi(x)^2)\zeta(x)^2\,dx\\
\nonumber
&\quad +C\int_{\Omega}|D\zeta(x)|^2\,dx+ \varepsilon C \int_{\Omega} |D\bu(x)|^2\zeta(x)^2\,dx
\end{align}
with constants depending on $n, \Lambda,\gamma, M,$ and $ \varepsilon.$
Summing up \eqref{Gradient} over $\alpha$ from $1$ to $N,$   fixing $\varepsilon$ small enough   and moving the last term to the left-hand side we obtain
\begin{align}\label{Gradient2}
\nonumber
\int_{\Omega} & |D\bu(x)|^2\zeta(x)^2\,dx  \leq C\int_{\Omega}(1+\psi(x)+\varphi(x)^2) \zeta(x)^2\, dx\\
&+C \sum_{\alpha=1}^N\int_{\Omega} |f^\alpha(x)| \zeta(x)^2\,dx
  +C\int_{\Omega}|D\zeta(x)|^2\,dx\,.
\end{align}
Then, by the definition of $\zeta$ and by \eqref{regularity} we have

\begin{align*}
& \int_{\B_{2\rho}} \big( 1+  \psi(x)+   \varphi(x)^2   \big)\,dx  \leq C  \big[\rho^n+ \rho^{n-\frac{n-\mu}{q}}\|\psi\|_{q,\mu;\Omega}\\
&\qquad\qquad   + \rho^{n-\frac{2(n-\lambda)}{p}}\|\varphi\|^2_{p,\lambda;\Omega}\big]\\
&\sum_{\alpha=1}^N \int_{\B_{2\rho}} |f^\alpha(x)| \,dx \leq 
C \rho^{n-\frac{n-\mu}{q}}\|\ff\|_{q,\mu;\Omega} \\
&\int_{\B_{2\rho}}|D\zeta(x)|^2\,dx\leq  C \rho^{n-2}\,.
\end{align*}
Hence
\begin{equation}\label{Gradient3}
\int_{\B_\rho}|D\bu |^2\,dx \leq C \rho^{\lambda_0}
\end{equation}
with $\lambda_0=\min\{n-2, n-\frac{2(n-\lambda)}{p}, n-\frac{n-\mu}{q} \}=n-2$ and the constant depends on known quantities.

\smallskip

Let $\B_{\rho}(x_0)\cap \partial\Omega\not=\emptyset.$ Then we extend $u^\alpha$ and the given functions $f^\alpha,\varphi,$ and $\psi$ as zero in $\Omega^c$  and consider the test functions 
$$
\phi^\alpha(x)=(e^{|u^\alpha(x)|}-1)\zeta(x)^2\sign u^\alpha(x)
$$ 
where $\zeta(x)$ is the cut-off function defined above. Thus \eqref{weak} gives
\begin{align*}
\sum_{i=1}^n& \int_{\Omega} A_i^\alpha(x,\bu,D\bu) e^{|u^\alpha(x)|} D_iu^\alpha(x)\zeta(x)^2\,dx\\
&\   = \int_{\Omega} f^\alpha(x)  (e^{|u^\alpha(x)|}-1)\zeta(x)^2 \sign u^\alpha(x)\,dx\\
&\quad - \sum_{i=1}^n
 \int_{\Omega} A_i^\alpha(x,\bu,D\bu) (e^{|u^\alpha(x)|}-1)2\zeta(x)D_i\zeta(x) \sign u^\alpha(x)\,dx\\
&\quad -  \int_{\Omega} b^\alpha(x,\bu,D\bu) (e^{|u^\alpha(x)|}-1)\zeta(x)^2\sign u^\alpha(x)\,dx\,.
\end{align*}
Hence the  conditions  \eqref{contr1} and \eqref{coercivity2} give
\begin{align*}
 \gamma \int_{\Omega}& |D u^\alpha(x)|^2\zeta(x)^2\,dx\leq 
\Lambda e^{M}\int_{\Omega} \varphi(x)^2 \zeta(x)^2 \,dx\\
&\  +e^{M}\int_\Omega|f^\alpha(x)|\zeta(x)^2\,dx\\
&\   +2 n\Lambda e^{M} \int_{\Omega} (\varphi(x)+ |\bu|^{\frac{n}{n-2}}+|D\bu|)\zeta(x) |D\zeta(x)|\,dx
\end{align*}
and to get the desired
 estimate \eqref{Gradient4} we argue as above.
\end{pf}

\end{document}